\documentclass[fullpage,twoside]{article}
\usepackage{fancyhdr,graphicx}
\usepackage{amsfonts}
\usepackage{booktabs}
\usepackage{array}
\usepackage{amsmath}
\usepackage{graphicx}
\usepackage{amsmath,bm}
\def \x{{\mbox{\boldmath $x$}}}

\def \hxi{{\mbox{\boldmath $\xi$}}}
\def \intd{{\rm d}}

\font\euler=eusm10

\def \M{\mbox{\euler M}}
\def \N{\mbox{\euler N}}

\newtheorem{definition}{Definition}[section]
\newtheorem{thm}{Theorem}[section]

\begin{document}
\title{\vspace{-1cm}  Uncertain Programming Model for Multi-item Solid Transportation Problem}
\author{Hasan Dalman$^{1,}$ \\
{\small\em $^{1,}$ Deparment of Mathematics Engineering, Yildiz Technical University, Esenler 34210, Turkey}\\
{\small\em hsandalman@gmail.com}}
\date{}
\maketitle
\maketitle \thispagestyle{fancy} \label{first}
\vspace{-5mm} \pagestyle{plain}
\newsavebox{\mytable}

\begin{abstract}

In this paper, an uncertain Multi-objective Multi-item Solid Transportation Problem (MMSTP) based on uncertainty theory is presented. In the model, transportation costs, supplies, demands and conveyances parameters are taken to be uncertain parameters. There are restrictions on some items and conveyances of the model. Therefore, some particular items cannot be transported by some exceptional conveyances.    Using the advantage of uncertainty theory, the MMSTP is first converted into an equivalent deterministic MMSTP. By applying convex combination method and  minimizing distance function method, the deterministic MMSTP is reduced into  single objective programming problems. Thus, both single objective programming problems are solved using Maple 18.02 optimization toolbox. Finally, a numerical example is given to illustrate the performance of the models.  


\end{abstract}

\section{Introduction}
The classical transportation problem (TP) is a simple mathematical programming problem in operational research, in which two types of constraints are considered i.e., source constraint and destination constraint. Nevertheless, in practice, we frequently are exposed to further restrictions besides of source constraint and destination constraint, such as goods constraint or transportation mode constraint. In such cases, the classical TP is transformed into a three-dimensional transportation problem, and it is called a Solid Transportation problem (STP). 

As a continuation of the classical TP, the STP is introduced by Haley \cite{1} in 1962. He is first proposed a solution procedure for a STP. 

Most of the real-world problems are generally characterized by multiple and conflicting criteria. Such conditions are usually formulated by optimizing multiple objective functions. In addition, the parameters of such models can be uncertain due to several uncontrollable factors. Therefore, the constructing of the real-world problems occur uncertainty with hesitation. Thus, we usually need to attach the uncertain parameters (fuzzy, interval or stochastic) to the models. Due to their facility to handle with the high level of uncertainty, fuzzy, interval and stochastic theory have been further used in various real-world applications, containing data mining \cite{2, 3, 4, 5, 6}, intelligent control \cite{7}, transportation \cite{8, 9, 10}, decision making \cite{11, 12}, so far.    

Until now, many researchers have investigated different types of Multi-objective Solid Transportation Problems (MSTPs) under uncertainty. Bit et al.\cite{13} proposed a fuzzy programming model for a MSTP. Jimenez and Verdegay \cite{14} formulated two types of the STP under uncertainty; that is, the supplies, demands and conveyance capacities in the model are considered interval and fuzzy numbers, respectively. In paper \cite{15}, an evolutionary algorithm based on the parametric approach is applied to obtain the optimal solutions of fuzzy STPs. In \cite{17}, a solution procedure  is introduced for solving the fixed charge MSTPs with type-2 fuzzy variables.

These days in a very frequently changing market, the trade of a single item does not pay much profit to a retail dealer. Consequently, all traders in the fields of transportation do the trade of several items in practice. Kundu et al. \cite{16} studied the fuzzy Multi-objective Multi-item Solid Transportation Problems (MMSTPs) with several items. They used the minimum fuzzy number to obtain the expected optimal solution of fuzzy MMSTPs. Dalman et al. \cite{18} proposed an interval fuzzy programming method to obtain the efficient solution of fuzzy MMSTPs. Furthermore, various types of MMSTPs under uncertainty are investigated by Das et al. \cite{19}, Giri et al. \cite{20}, Yang et al. \cite{21}.    

According to Liu et al. \cite{22}, when we use probability theory, enough historical information is required to forecast the probability distribution. In various conditions, there are no examples available to forecast the probability distribution and then researchers are invited to estimate the degree of belief that each condition may occur. To deal with the degree of belief, uncertainty theory introduced by Liu \cite{23} in 2007 and improved by Liu \cite{24} in 2010. Up to now, uncertainty theory has been applied to many other fields. On the application of uncertainty theory, the interested readers may consult Liu \cite{22, 24, 25, 26}, Wang et al. \cite{27} and Gao \cite{28}. 

Uncertain programming is a kind of mathematical programming including uncertain parameters. The theory of uncertain programming is presented by Liu \cite{29} in 2009. It is a practical tool handling decision process, including the degree of belief. Up to now, uncertain programming used to optimize machine scheduling, vehicles routing, and project scheduling problems \cite{33}. Liu and Yao \cite{30} suggested an uncertain goal programming method for solving the uncertain multi-objective programming problem. Liu and Yao \cite{31} presented some programming method to solve uncertain multilevel programming problems. Zhou et al. \cite{32} suggested the compromise programming model for solving uncertain multi-objective programming problems. Zhong et al. \cite {33} suggested an interactive satisfied method for an uncertain multi-objective programming problem. Cui and Sheng \cite{34} formulated an expected constrained programming model for the STP based on uncertainty theory.   

As far as I know, no work has been studied for the MMSTP based on uncertainty theory. Thus, this paper presents a MMSTP with uncertain parameters. Using the expected value of each objective function under the chance constraints, the model is transformed into a deterministic MMSTP. With the use of convex combination method and distance minimizing method, the deterministic MMSTP is reduced to single objective programming problems. Then single objective programming problems are solved to obtain the optimal solutions of the considered problem by using Maple 18.02 optimization toolbox. Finally, a numerical example is presented to illustrate the solution procedures.     

This paper is constructed as follows: Section 2 presents some definitions and theorems about uncertainty theory. Section 3 presents uncertain programming and its solution procedure. In Section 4, a MMSTP based on uncertainty theory is formulated and then the solution procedure is illustrated by suitable numerical example in Section 5.  

\section{Preliminary}

In this section, we will give some basic definitions and notations about uncertainty theory.

\begin{definition} (Liu {\cite{23}})
Let  $\mbox{\euler L}$ be a $\sigma$-algebra on a nonempty set  $\Gamma$. A set function $\M$ is called an \emph{uncertain measure} if it satisfies the following axioms:

\vskip 0.15mm\noindent {Axiom 1.}~(\emph{Normality Axiom}) $\M\{\Gamma\}=1$;
\vskip 0.15cm\noindent {Axiom
2.}~(\emph{Duality Axiom}) $\M\{\Lambda\}+\M\{\Lambda^{c}\}=1$ for any
 $\Lambda \in \mbox{\euler L} $;
\vskip 0.15cm\noindent {Axiom 3.}~(\emph{Subadditivity
Axiom}) For every countable sequence of
$\{\Lambda_i\}\in \mbox{\euler L}$, we have
\begin{equation*}\displaystyle
\M\left\{\bigcup_{i=1}^{\infty}\Lambda_i\right\}\le\sum_{i=1}^{\infty}\M\{\Lambda_i\}.
\end{equation*}

The triplet $(\Gamma,\mbox{\euler L},\M)$ is called an {\em uncertainty
space}, and each element $\Lambda$ in $\mbox{\euler L} $ is called an {\em event}. In addition, in order to obtain an
uncertain measure of compound event, a product uncertain measure is defined
by Liu \cite{35} by the following product axiom:

\vskip 0.15cm\noindent{Axiom 4.}~({\em Product Axiom}) Let $(\Gamma_k,\mbox{\euler
L}_k,\M_k)$ be uncertainty spaces for $k=1,2,\cdots$ The product
uncertain measure $\M$ is an uncertain measure satisfying

\begin{equation*} \M\left\{\prod_{k=1}^{\infty}\Lambda_k\right\}=\bigwedge_{k=1}^{\infty}\M_k\{\Lambda_k\}\end{equation*}
where $\Lambda_k$ are
arbitrarily chosen events from $\mbox{\euler L}_k$ for
$k=1,2,\cdots$, respectively.
\end{definition}

\begin{definition} (Liu {\cite{23}}) An {\em uncertain variable} $\xi$ is a measurable function from
an uncertainty space $(\Gamma,\mbox{\euler L},\M)$ to the set of
real numbers, i.e., for any Borel set $B$ of real numbers, the set
 \begin{equation*}
\{\xi\in B\}=\{\gamma\in \Gamma|\xi(\gamma)\in B\}
\end{equation*}
is an event.
\end{definition}

\begin{definition} (Liu {\cite{23}})
The {\em uncertainty distribution} $\Phi$ of an uncertain variable $\xi$
is defined by
\begin{equation*} \Phi(x)=\M\left\{\xi\le
x\right\},\quad\forall x\in\Re.
\end{equation*}
\end{definition}

\begin{definition}(Liu {\cite{23}})
Let $\xi$ be an uncertain variable. The {\em expected value} of $\xi$  is defined by
\begin{equation*}
E[\xi]=\int_{0}^{+\infty}\M\{\xi\geq r\}\intd r-\int_{-\infty}^{0}\M\{\xi\leq r\}\intd r
\end{equation*}
{\em provided that at least one of the above two integrals is finite.} 
An uncertain variable $ \xi $ is called linear if it has a linear uncertainty distribution

\[\Phi \left( x \right) = \left\{ {\begin{array}{*{20}{l}}
	{0,\,\,\,\,\,\,\,\,\,\,\,\,\,\,\,\,\,\,\,\,\,\,\,\,\,\,\,\,\,\,\,\,\,\,\,\,\,\,\,\,\,\,\,\,x \le a}\\
	{\left( {x - a} \right)/\left( {b - a} \right),\,\,\,\,a \le x \le b}\\
	{1,\,\,\,\,\,\,\,\,\,\,\,\,\,\,\,\,\,\,\,\,\,\,\,\,\,\,\,\,\,\,\,\,\,\,\,\,\,\,\,\,\,\,\,\,x \ge b}
	\end{array}} \right.\]
denoted by $ \L\left( {a,b} \right) $   where $ a $   and  $ b $   are real numbers with $a < b.$ Suppose that $ \xi_1 $ and $ \xi_2 $  are independent linear uncertain variables $ \L\left( {a_1,b_1} \right) $ and $ \L\left( {a_2,b_2} \right) .$  Then the sum $ {\xi _1} + {\xi _2} $  is also a linear uncertain variable $ \L\left( {{a_1} + {a_2},{b_1} + {b_2}} \right). $
\end{definition}
\begin{definition} (Liu {\cite{24}})
	An uncertainty distribution  $\Phi(x)$ is said to be {\em regular }if it is a continuous and
	strictly increasing function with respect to $x$ at which $0 < \Phi(x) < 1$, and
	\begin{equation*}
	\lim_{x\rightarrow -\infty} \Phi(x)=0,\ \ \  \  \lim_{x\rightarrow +\infty} \Phi(x)=1.
	\end{equation*}
	
\end{definition}

\begin{definition}
	Let $\xi$ be an uncertain variable with a regular uncertainty distribution $\Phi(x)$. If the expected value is available, then
	\[E[\xi]=\int_{0}^{1}{\Phi^{-1}(\alpha)}\intd \alpha \]
	where $\Phi^{-1}(\alpha)$ is the {\em inverse uncertainty distribution} of $\xi$.

An uncertain variable $ \xi $  is called normal if it has a normal uncertainty distribution 
\[ \Phi \left( x \right) = \left( {1 + \exp \left( {\pi \frac{{\left( {e - x} \right)}}{{\sqrt 3 \sigma }}} \right)} \right), - \infty  < x < \infty ,\sigma  > 0 \]
denoted by $ {\rm N}\left( {e,\sigma } \right). $Assume that $ \xi_1 $and $ \xi_2 $  are independent normal uncertain variables $ {\rm N}\left( {e_1,\sigma_1 }  \right)$   and $ {\rm N}\left( {e_2,\sigma_2 }  \right)$  Then the sum $ \xi_1+\xi_2 $ is also a normal uncertain variable ${\rm N}\left( {{e_1} + {e_2},{\sigma _1} + {\sigma _2}} \right).$

Note that the expected value of linear uncertain variable $ L\left( {a,b} \right) $   and ${\rm N}\left( {e,\sigma } \right)$ normal uncertain variable  are $ \frac{{a + b}}{2} $ and $ e $  respectively.

In addition, the inverse uncertainty distribution of a normal uncertain variable ${\rm N}\left( {e,\sigma } \right)$ is 
\[{\Phi ^{ - 1}}\left( x  \right) = {e} + \frac{{{\sigma}\sqrt 3 }}{\pi }\ln \frac{{{x}}}{{1 - {x}}}\]]
\end{definition}
\begin{thm} Assume $\xi_1,\xi_2,\cdots,\xi_n$ are independent uncertain variables
with regular uncertainty distributions $\Phi_1,\Phi_2,\cdots,\Phi_n$,
respectively. If the function $f(x_1,x_2,\cdots,x_n)$ is strictly
increasing with respect to $x_1,x_2,\cdots$, $x_m$ and strictly
decreasing with respect to $x_{m+1},x_{m+2},\cdots,x_n$, then
$\xi=f(\xi_1,\xi_2,\cdots,\xi_n)$ has an
inverse uncertainty distribution
\begin{equation*}
\Psi^{-1}(\alpha)=f\left(\Phi_1^{-1}(\alpha),\cdots,
\Phi_m^{-1}(\alpha),\Phi_{m+1}^{-1}(1-\alpha),\cdots,\Phi_n^{-1}(1-\alpha)\right).
\end{equation*}
In addition, Liu and Ha \cite{36} proved that the uncertain variable $\xi$ has an expected value
\begin{equation*}
E[\xi]=\int_{0}^{1}f\left(\Phi_1^{-1}(\alpha),\cdots,
\Phi_m^{-1}(\alpha),\Phi_{m+1}^{-1}(1-\alpha),\cdots,\Phi_n^{-1}(1-\alpha)\right)\intd\alpha.
\end{equation*}
\end{thm}

\section{Uncertain Programming}

Assume that $x = \left( {x_1^{},x_2^{},...,x_n^{}} \right)$ is $ n- $ dimensioal decision vector, $\xi  = \left( {\xi _1^{},\xi _2^{},...,\xi _n^{}} \right)$ is $ n- $ dimensional uncertain vector, $f(\x,\hxi)$ and ${g_j}(\x,\hxi)\le0$ are the objective and constraint functions, respectively. However, these functions do not specify deterministic functions. It is inherently expected that the uncertain  constraints hold with defined confidence levels $\alpha  = \left( {\alpha _1^{},\alpha _2^{},...,\alpha _p^{}} \right).$ Then we have a set of distinct chance constraints,
{\begin{equation*} \M\left\{ {{g_j}\left( {x,\xi } \right) \le 0} \right\} \ge {\alpha _j}\end{equation*}}
Here, function $f(\x,\hxi)$ cannot be optimized naturally. So we optimize its expected value, i.e.,
{\begin{equation*}\min_{\x}E[f(\x,\hxi)].\end{equation*}}
In order to optimize an uncertain programming model, consider the following programming model which is introduced by  Liu \cite{29};
{\begin{equation*}
\left\{\begin{array}{l}
\min\limits_{\x} E[f(\x,\hxi)] \\[0.2cm]
\mbox{subject to:} \\[0.1cm]
\qquad \M\{g_j(\x,\hxi)\le 0\}\ge\alpha_j,\quad j=1,2,\cdots,p
\end{array}\right.
\end{equation*}}
where $\x$ is a decision vector, $\hxi$ is an uncertain vector, $f$
is an objective function, and $g_j$ are constraint functions for $j=1,2,\cdots,p$.

In this uncertain programming model, (From theorem 2.1) we can determine the expected value  {$E[f(\x,\xi_1,\xi_2,\cdots,\xi_n)]$} of uncetain objective function as follows: 
\begin{equation*}
\int_0^1f(\x,\Phi_1^{-1}(\alpha),\cdots,\Phi_m^{-1}(\alpha),
\Phi_{m+1}^{-1}(1-\alpha),\cdots,\Phi_n^{-1}(1-\alpha))\intd\alpha.
\end{equation*}
The inverse uncertainty distribution of
$f(\x,\xi_1,\xi_2,\cdots,\xi_n)$ is
\begin{equation*}\Psi^{-1}(\x,\alpha)=f(\x,\Phi_1^{-1}(\alpha),\cdots,\Phi_m^{-1}(\alpha),
\Phi_{m+1}^{-1}(1-\alpha),\cdots,\Phi_n^{-1}(1-\alpha))
\end{equation*}

In a similar manner, the chance constraint is:
{\begin{equation*}
	\M\left\{g(\x,\xi_1,\xi_2,\cdots,\xi_n)\le
	0\right\}\ge\alpha\end{equation*}} holds if and only if
{\begin{equation*}
	g(\x,\Phi_1^{-1}(\alpha),\cdots,\Phi_k^{-1}(\alpha),
	\Phi_{k+1}^{-1}(1-\alpha),\cdots,\Phi_n^{-1}(1-\alpha))\le 0.
	\end{equation*}}
The inverse uncertainty distribution of
$g(\x,\xi_1,\xi_2,\cdots,\xi_n)$ is
$$\Psi^{-1}(\x,\alpha)=g(\x,\Phi_1^{-1}(\alpha),\cdots,\Phi_k^{-1}(\alpha),
\Phi_{k+1}^{-1}(1-\alpha),\cdots,\Phi_n^{-1}(1-\alpha)).$$

In actual decision making conditions, a major concern is that
most decision problems involve multiple objectives. In order to construct uncertain programming models with multiple objectives, we extend the single objective programming to obtain the following uncertain multi-objective programming model.
\begin{equation}
	\left\{\begin{array}{l} \min\limits_{{\mbox{\footnotesize\boldmath
					$x$}}}
		\left(E[f_1(\x,\hxi)],E[f_2(\x,\hxi)],\cdots,E[f_q(\x,\hxi)]\right) \\[0.1cm]
		\mbox{subject to:} \\[0.1cm]
		\qquad \M\{g_j(\x,\hxi)\le 0\}\ge\alpha_j,\quad j=1,2,\cdots,p
	\end{array}\right.
\end{equation}\\ 
where $E[f_1(\x,\hxi)],E[f_2(\x,\hxi)],\cdots,E[f_q(\x,\hxi)]$
are real-valued objective functions.

\begin{definition}
An element $\x^*$ is called a feasible decision to the above uncertain multi-objective programming if $\qquad \M\{g_j(\x,\hxi)\le 0\}\ge\alpha_j,\quad$ for $ j=1,2,...,p.$
When the objectives are in conflict, there is no optimal solution $\x^*$ that simultaneously minimizes all the objective functions $ E[f_i(\x^*,\hxi)],$$ i=1,2,...,q.$ Therefore, the idea of Pareto solution is defined.
\end{definition}

\begin{definition}
A feasible decision $\x^*$ is called a Pareto optimal solution of the uncertain multi-objective programming model if there is no feasible decision $ x $ such that 
{$E[f_i(\x,\hxi)]\le E[f_i(\x*,\hxi)]$} for all $ i=1,2,...,q.$ and {$E[f_k(\x,\hxi)]<E[f_k(\x*,\hxi)]$} for at least one index $ k.$

\end{definition}

In order to optimize the uncertain programming models, various methods have been suggested. In this paper, two existing compromise programming method will employ to achieve the Pareto optimal solution of uncertain multi-objective programming problems.   

\subsection{The  Convex Combination Method}
Consider the above uncertain multi-objective programming model (1). The compromise programming concept is employed to achieve Pareto solutions in model (1).    

By weighting the objective functions, we have the following compromise programming model.

\begin{equation}
\left\{\begin{array}{l}
\min\limits_{{\mbox{\footnotesize\boldmath $x$}}}\sum\limits_{i=1}^m\lambda_iE[f_i(\x,\hxi)]\\[0.3cm]
\mbox{subject to:}\\[0.1cm]
\qquad \M\{g_j(\x,\hxi)\le 0\}\ge\alpha_j,\quad j=1,2,\cdots,p
\end{array}\right.
\end{equation}
where the weights $\lambda_1,\lambda_2,...,\lambda_q$ are
nonnegative numbers with $\lambda_1+\lambda_2+...+\lambda_q=1.$

The above method changes multiple objectives into a summed objective function by multiplying each objective function by a weighting factor and gathering all weighted objective functions.
\begin{thm} 
Let $ x* $ be an optimal solution of the convex combination model (2).Then $ x* $ must be a Pareto optimal solution of the uncertain multi-objective programming model (1).
\end{thm}

Assume that the optimal solution $ x* $ is not a Pareto optimal solution of the uncertain multi-objective programming model (1). So, there must exist a feasible solution $ x $ such that {$E[f_i(\x,\hxi)]\le E[f_i(\x*,\hxi)]$} for all $ i=1,2,...,q.$ and {$E[f_k(\x,\hxi)]<E[f_k(\x*,\hxi)]$} for at least one index $ k=$
Therefore, we write for $ \lambda>0,$  {$\lambda_iE[f_i(\x,\hxi)]\le \lambda_iE[f_i(\x*,\hxi)]$} for all $ i=1,2,...,q.$ and {$ \lambda_kE[f_k(\x,\hxi)]<\lambda_kE[f_k(\x*,\hxi)] $} for at least one index $ k$.

As a result, we obtain
\[ \sum\limits_{i = 1}^q \lambda _i E[f_i(\x,\hxi)]\le \sum\limits_{i = 1}^q \lambda _i E[f_i(\x*,\hxi)] \]
which implies that $ x* $ is not an optimal solution of model (2). So, $ x* $ is a Pareto optimal solution of the uncertain programming model (1).

\subsection{Minimizing Distance Function}
The second way is related to minimizing the distance function from a solution   
\[ (E[f_1(\x,\hxi)],E[f_2(\x,\hxi)],\cdots,E[f_q(\x,\hxi)]) \]
to an ideal vector $ (E_1^*,E_2^*,...,E_q^*) $ where $ E_i^* $ are the optimal values of the $i-$th objective functions without considering other objectives, $i=1,2,...,q$, respectively i.e,
\begin{equation}
\left\{\begin{array}{l} \min\limits_{{\mbox{\footnotesize\boldmath
			$x$}}}
\left(\sqrt{(E[f_1(\x,\hxi)]-E_1^*)^2+(E[f_2(\x,\hxi)]-E_2^*)^2+,\cdots,+(E[f_q(\x,\hxi)]-E_q^*)^2}\right) \\[0.1cm]
\mbox{subject to:} \\[0.1cm]
\qquad \M\{g_j(\x,\hxi)\le 0\}\ge\alpha_j,\quad j=1,2,\cdots,p
\end{array}\right.
\end{equation}\\
\begin{thm} 
Let $ x* $ be an optimal solution of the uncertain programming model (3).Then $ x* $ must be a Pareto optimal solution of the uncertain multi-objective  programming model (1).
\end{thm}

Assume that the optimal solution $ x* $ is not a Pareto optimal solution of the uncertain programming model (1). So, there must exist a feasible solution $ x $ such that {$E[f_i(\x,\hxi)]\le E[f_i(\x*,\hxi)]$} for all $ i=1,2,...,q.$ and {$E[f_k(\x,\hxi)]<E[f_k(\x*,\hxi)]$} for at least one index $ k$.
Since  $ E_i^*$ are the optimal solution of the $i$th objective function without considering other objective functions $ i=1,2,...,q.$, we obtain  {$ E_k^*\le E[f_k(\x,\hxi)]\le E[f_k(\x*,\hxi)]$}  and {$ E_i^*\le E[f_i(\x,\hxi)]\le E[f_i(\x*,\hxi)] $}.
Therefore, we write $ (E[f_1(\x,\hxi)]-E_1^*)^2+...+(E[f_k(\x,\hxi)]-E_k^*)^2+...+(E[f_q(\x,\hxi)]-E_q^*)^2\le (E[f_1(\x*,\hxi)]-E_1^*)^2+...+(E[f_k(\x*,\hxi)]-E_k^*)^2+...+(E[f_q(\x*,\hxi)]-E_q^*)^2 $.

Hereby,  $ x* $ is not an optimal solution of  model (3). A contradiction shows that $ x* $ is a Pareto optimal solution of  model (1).

\section{Uncertain Programming Model For Multi-objective Multi-item Solid Transportation Problem}

In this section, a MMSTP with uncertain parameters will be formulated. In this model, following notations are used.

$ i $, $ {1, 2, . . . , m} $  is the index of sources,

$ j $, $1, 2, . . . , n $ is the index of destinations,

$ k $, $ {1, 2, . . .,l} $ is the index for conveyances,

$ p $, $ {1, 2, . . .,r} $ is the index of items,

$ t $,  $ {1, 2, . . .,K} $ is the index of objectives,

$ a^p_i $ is the amount of products of item $ p $ at source,

$ b^p_j $ is the demand of products of item $ p $ at destination $ j $,

$ e_k $ is the total transportation capacity of conveyance $ k $,

$c^{tp}_{ijk} $ is the cost for transporting one unit of item  $ p $ from source $ i $ to destination $ j $ by conveyance $ k $,

$x^{p}_{ijk} $  is the amount of item $ p $ to be transported from source $ i $ to destination $ j $ with the aid of conveyance $ k $.

\subsection{Transportation cost}

A MMSTP with $p$ objectives in which  $r$ distinct  items are carried from $ m$ sources to $ n $ destinations by $ l $ different conveyances. 

If the transportation activity occurs, the transportation cost will be paid proportionally. In such case, the total transportation cost is contsructed as
\[{f_t(x)} = \sum\limits_{p = 1}^r\sum\limits_{i = 1}^m {\sum\limits_{j = 1}^n {\sum\limits_{k = 1}^l {{c^{tp}_{ijk}}{x^p_{ijk}}} }, t=(1,2,...,K) } \]. 
\subsection{Model Constraints}

The first constraint of the model is the total quantity carried from source $ i $ is no more than $ a^{p}_i .$ Then, we obtain
\[  \sum\limits_{j = 1}^n{\sum\limits_{k = 1}^l} {{x^{p}_{ijk}} \le a^{p}_i,{i = 1,2,...,m};{p = 1,2,...,r}} \]

The second constraint of the model is the total quantity carried from sources should meet the demand of destination $ j. $ Then we obtain;
\[  \sum\limits_{i = 1}^m{\sum\limits_{k = 1}^l} {{x^{p}_{ijk}} \ge b^{p}_j,{j = 1,2,...,n};{p = 1,2,...,r}} \]

The third constraint of the model is the total quantity carried by conveyance k is not more than its transportation capacity. Thus, we obtain
\[  \sum\limits_{p = 1}^r{\sum\limits_{i = 1}^m}{\sum\limits_{k = 1}^l} {{x^{p}_{ijk}} \le e_k,{k = 1,2,...,l}} \]

Using the above notations, the MMSTP can be formulated as follows:
1 \begin{equation}
\left\{ \begin{array}{l}
{f_t\left({x}\right)} = \sum\limits_{p = 1}^r\sum\limits_{i = 1}^m\sum\limits_{j = 1}^n\sum\limits_{k = 1}^l {c_{ijk}^{tp}}{x_{ijk}^{p}}, t=(1,2,...,K)\\
s.t.\left\{ \begin{array}{l}
\sum\limits_{j = 1}^n{\sum\limits_{k = 1}^l} {{x^{p}_{ijk}} \le a^{p}_i,{i = 1,2,...,m};{p = 1,2,...,r}} \\
\sum\limits_{i = 1}^m{\sum\limits_{k = 1}^l} {{x^{p}_{ijk}} \ge b^{p}_j,{j = 1,2,...,n};{p = 1,2,...,r}} \\
\sum\limits_{p = 1}^r{\sum\limits_{i = 1}^m}{\sum\limits_{k = 1}^l} {{x^{p}_{ijk}} \le e_k,{k = 1,2,...,l}} \\
x^{p}_{ijk} \ge 0,  \forall pijk
\end{array} \right.
\end{array} \right.
\end{equation}

The above model (4) is formulated under deterministic conditions; that is, the variables of the above model are all deterministic quantities. Due to the deficiency of information, we conventionally encounter with the uncertain event in building the mathematical model. Therefore, we generally add the uncertain parameters to the model.    

In this paper, we suppose that the unit cost, the capacity of each source and that of each destination are all uncertain variables. These variables are $\xi_{ijk}^{tp},\tilde a_{i}^{p},\tilde b_{j}^{p},\tilde e_k$, respectively.

Therefore, the above model (4) is converted to its uncertain programming model. It is called an expected constrained programming model. Here, the essential scheme of an uncertain transportation model is to optimize the expected value of objective function under the chance constraints. 
So, the expected value programming model under chance constraints is formulated as:  

\begin{equation}
\left\{ \begin{array}{l}
E[f_t(\x,\hxi)]=\min {\rm E}\left[{\sum\limits_{p = 1}^r {\sum\limits_{i = 1}^m {\sum\limits_{j = 1}^n {\sum\limits_{k = 1}^l {\left( {\xi _{ijk}^{tp}x_{ijk}^p} \right)} } } } } \right], t=(1,2,...,K)\\
s.t.\left\{ {\begin{array}{*{20}{l}}
{{\rm M}\left\{ \sum\limits_{j = 1}^n\sum\limits_{k = 1}^lx_{ijk}^{p}  -{\tilde a}_i^{p} \le 0 \right\} \ge {\gamma _i^{p}},i = 1,2,...,m}\\
{{\rm M}\left\{ {{{\tilde b}_j} - \sum\limits_{k = 1}^m \sum\limits_{k = 1}^l {{x_{ijk}^{p}} }  \le 0} \right\} \ge {\beta _j^{p}},j = 1,2,...,n.}\\
{{\rm M}\left\{ \sum\limits_{p = 1}^r\sum\limits_{i = 1}^m \sum\limits_{j = 1}^n {x_{ijk}^{p} - {{\tilde e}_k} \le 0} \right\} \ge {\delta _k},k = 1,2,...,l.}
\end{array}} \right.
\end{array} \right.
\end{equation}
where ${\alpha _i^{p}}$,${\beta _j^{p}}$,${\delta _k}$ are specified confidence levels of constraints for $i = 1,2,...,m$,$j = 1,2,...,n$,$k = 1,2,...,l, p=1,2,...,r$. 

In the above model, the following conditions are satisfied;

The first constraint implies that total amount transported from source $ i $ should be no more than its supply capacity at the confidence level  ${\alpha _i}$.

The second constraint implies that the total amount transported from source $ i $ should satisfy the requirement of destination $ j $ at the credibility level ${\beta _j}$.

The third constraint states that the total amount transported by conveyance $ k $ should be no more than its transportation capacity at the confidence level ${\delta _k}$.

\subsection{Deterministic Equivalences of Models}
In this section, we will show the deterministic equivalents of model (5) under some specific situations.
\begin{thm}
Assume that ${\xi^{tp}_{ijk}},{\tilde a_i^{p}},{\tilde b_j^{p}},{\tilde e_k}$ are independent uncertain variables with uncertainty distributions ${\Phi _\xi }_{ijk},{\Phi _{{{\tilde a}_i}}},{\Phi _{{{\tilde b}_j}}},{\Phi _{{{\tilde e}_k}}}.$ Then the above model (5) is converted into the following  equivalent model:
\end{thm}
\begin{equation}
\begin{array}{*{20}{l}}
{\min \sum\limits_{p = 1}^r \sum\limits_{i = 1}^m {\sum\limits_{j = 1}^n {\sum\limits_{k = 1}^l {{x_{ijk}}} } } \int\limits_0^1 {\Phi _{{\xi^{tp} _{ijk}}}^{ - 1}} \left( \alpha  \right)d\alpha}\\
{s.t.\left\{ {\begin{array}{*{20}{l}}
{\sum\limits_{j = 1}^n {\sum\limits_{k = 1}^l {{x^{p}_{ijk}}} }  - \int\limits_0^1 {\Phi _{{{\tilde \gamma^{p} }_i}}^{ - 1}} (1 - {\gamma^{p} _i})\left( \alpha  \right)d\alpha  \le 0,i = 1,2,...,m; p=1,2,...,r}\\
{\int\limits_0^1 {\Phi _{{{\tilde \beta^{p} }_j}}^{ - 1}} \left( \alpha  \right)d\alpha  - \sum\limits_{k = 1}^m {\sum\limits_{k = 1}^l {{x^{p}_{ijk}}} }  \le 0,j = 1,2,...,n;  p=1,2,...,r}\\
{\sum\limits_{p = 1}^r\sum\limits_{i = 1}^m {\sum\limits_{j = 1}^n {{x^{p}_{ijk}}} }  - \int\limits_0^1 {\Phi _{{{\tilde \delta }_k}}^{ - 1}} \left( \alpha  \right)d\alpha  \le 0,k = 1,2,...,l}
\end{array}} \right.}
\end{array}\\
\end{equation}
\begin{proof}
Since ${\xi^{tp} _{ijk}},{\tilde a^{p}_i},{\tilde b^{p}_j},{\tilde e_k}$ are independent uncertain variables with uncertainty distributions ${\Phi _{{\xi^{tp} _{ijk}}}},{\Phi _{{{\tilde a^{p}}_i}}},{\Phi _{{{\tilde b^{p}}_j}}},{\Phi _{{{\tilde e}_k}}}$, respectively. 

From the linearity of the expected value operator, we obtain   
\begin{equation}
{\rm E}\left[ {\sum\limits_{i = p}^r\sum\limits_{i = 1}^m {\sum\limits_{j = 1}^n {\sum\limits_{k = 1}^l {{\xi^{tp} _{ijk}}} {x^{p}_{ijk}}} } } \right] = \sum\limits_{i = p}^r\sum\limits_{i = 1}^m {\sum\limits_{j = 1}^n {\sum\limits_{k = 1}^l {{x^{p}_{ijk}}} } } {\rm E}\left[ {{\xi^{tp} _{ijk}}}  \right],t=1,2,...,K.
\end{equation}
where $E\left[{{\xi^{tp} _{ijk}}}\right]=\int\limits_0^1 {\Phi _{{\xi^{tp} _{ijk}}}^{ - 1}} \left( \alpha  \right)d\alpha$, $i = 1,2,...,m$,$j = 1,2,...,n$,$k = 1,2,...,l, t=1,2,...,K$. 
According to Theorems 2.1, we have:\\
\begin{equation}
\begin{array}{l}
{\rm M}\left\{ {\sum\limits_{j = 1}^n {\sum\limits_{k = 1}^l {{x^{p}_{ijk}}} }  - {{\tilde a^{p}}_i} \le 0} \right\} \ge {\gamma _i^{p}} \Leftrightarrow \sum\limits_{j = 1}^n {\sum\limits_{k = 1}^l {{x^{p}_{ijk}}} }  - \int\limits_0^1 {\Phi _{{{\tilde \gamma^{p}}_i}}^{ - 1}} (1 - {\gamma^{p} _i})\left( \alpha  \right)d\alpha  \le 0,i = 1,2,...,m; p=1,2,...,r\\
{\rm M}\left\{ {{{\tilde b^{p}}_j} - \sum\limits_{k = 1}^m {\sum\limits_{k = 1}^l {{x^{p}_{ijk}}} }  \le 0} \right\} \ge {\beta^{p} _j} \Leftrightarrow \int\limits_0^1 {\Phi _{{{\tilde \beta^{p} }_j}}^{ - 1}} \left( \alpha  \right)d\alpha  - \sum\limits_{k = 1}^m {\sum\limits_{k = 1}^l {{x_{ijk}}} }  \le 0,j = 1,2,...,n; p=1,2,...,r\\
{\rm M}\left\{ {\sum\limits_{p =1}^r\sum\limits_{i = 1}^m {\sum\limits_{j = 1}^n {{x^{p}_{ijk}}} }  - {{\tilde e}_k} \le 0} \right\} \ge {\delta _k} \Leftrightarrow \sum\limits_{p =1}^r\sum\limits_{i = 1}^m {\sum\limits_{j = 1}^n {{x^{p}_{ijk}}} }  - \int\limits_0^1 {\Phi _{{{\tilde \delta }_k}}^{ - 1}} (1 - {\delta_i})d\alpha  \le 0,\ k = 1,2,...,l.
\end{array}
\end{equation}

\end{proof}

Thus, model (4) is converted into its deterministic model (6). Then we can solve easily it to obtain the Pareto optimal solutions.

\section{A Numerical Example}
Suppose the multi-objective multi-item STP is to be carried by two distinct conveyances from four sources to three destinations and consisting of two different objective functions. Solve the problems to obtain the amount of goods to be shipped from source(s) to destination(s) so that the total demand at all the destinations is met at the minimum total cost and in the potential minimum time.   

In the following problem, all uncertain parameters are considered as normal uncertain variables. Also these parameters can be chosen as linear and/or zigzag uncertain parameters by the decision maker.

\begin{equation}
\begin{array}{l}
{\xi^{tp} _{ijk}}\sim \N\left( {e^{tp}_{ijk}},\sigma^{tp} _{ijk} \right)$, i = 1, 2,3,4, j = 1, 2,3,4,5,6, k = 1, 2,p = 1, 2, t=1,2.$\\
{{\tilde a^p}_i}\sim{\rm N}\left( {{e^p_i},{\sigma^p_i}} \right),i = 1,2,3, p=1,2.\\
{{\tilde b^p}_j}\sim{\rm N}\left( {{{e'^p}_j},{{\sigma '^p}_j}} \right),j = 1,2,3,4, p=1,2.\\
{{\tilde e}_k}\sim{\rm N}\left( {{{e''}_k},{{\sigma ''}_k}} \right),k = 1,2.
\end{array}
\end{equation}Data for the transportation problem is given in Tables 5.1, 5.2, 5.3, 5.4, 5.5, 5.6 and 5.7.
\begin{center}
 Table 5.1 :The transportation cost $ {\xi^{11}_{ijk}}$  for item 1 in the first objective\\
\begin{tabular}{cccccccccc} 
\hline $ {\xi^{11}_{ij1}}$&1&2&3&4&$ {\xi^{11}_{ij2}}$&1&2&3&4\\
\hline 
1&(10,2)&(9,1.5)&(12,2)&(8,1.5)&1&(7,2)&(5,1.5)&(5,2)&(6,1.5)\\
2&(8,1)&(9,1.5)&(11,2)&(10,1)&2&(5,1)&(6,1.5)&(4,2)&(6,1)\\
3&(8,1.5)&(17,1.5)&(6,1.5)&(10,1.5)&3&(4,1.5)&(7,1.5)&(6,1.5)&(5,1.5)\\
\hline
\end{tabular}\\
\end{center}
\begin{center}
	Table 5.2 :The transportation cost $ {\xi^{12}_{ijk}}$ for item 2 in the first objective\\
\begin{tabular}{cccccccccc} 
	\hline $ {\xi^{12}_{ij1}}$&1&2&3&4&$ {\xi^{12}_{ij2}}$&1&2&3&4\\
	\hline 
	1&(9,2)&(6,1.5)&(3,2)&(7,1.5)&1&(5,2)&(6,1.5)&(7,2)&(9,1.5)\\
	2&(9,1)&(10,1.5)&(11,2)&(10,1)&2&(5,1)&(5,1.5)&(3,2)&(3,1)\\
	3&(6,1.5)&(18,1.5)&(8,1.5)&(12,1.5)&3&(2,1.5)&(8,1.5)&(7,1.5)&(3,1.5)\\
	\hline
\end{tabular}\\
\end{center}
\begin{center}
	Table 5.3 :The transportation cost $ {\xi^{21}_{ijk}}$ for item 1 in the second objective\\
\begin{tabular}{cccccccccc} 
		\hline $ {\xi^{21}_{ij1}}$&1&2&3&4&$ {\xi^{21}_{ij2}}$&1&2&3&4\\
		\hline
1&(20,2)&(19,1.5)&(22,2)&(21,1.5)&1&(30,2)&(25,1.5)&(27,2)&(29,1.5)\\
	2&(23,1)&(22,1.5)&(18,2)&(17,1)&2&(26,1)&(26,1.5)&(28,2)&(32,1)\\
	3&(16,1.5)&(17,1.5)&(20,1.5)&(18,1.5)&3&(33,1.5)&(29,1.5)&(35,1.5)&(30,1.5)\\
	\hline
	\end{tabular}\\
\end{center}
\begin{center}
	Table 5.4 :The transportation cost $ {\xi^{22}_{ijk}}$ for item 2 in the second objective\\
	\begin{tabular}{cccccccccc} 
		\hline $ {\xi^{22}_{ij1}}$&1&2&3&4&$ {\xi^{22}_{ij2}}$&1&2&3&4\\
		\hline 
		1&(18,2)&(16,1.5)&(23,2)&(17,1.5)&1&(15,2)&(16,1.5)&(17,2)&(19,1.5)\\
		2&(19,1)&(20,1.5)&(21,2)&(20,1)&2&(24,1)&(25,1.5)&(23,2)&(23,1)\\
		3&(16,1.5)&(18,1.5)&(18,1.5)&(12,1.5)&3&(22,1.5)&(28,1.5)&(27,1.5)&(23,1.5)\\
		\hline
	\end{tabular}\\
\end{center}
\begin{center}
Table 5.5 The sources ${\tilde a^p_i}$\\
\begin{tabular}{cccccccc}
\hline $i$&1&2&3&&1&2&3\\
\hline 
${\tilde a^1_i}$&(32,1.5)&(35,1.5)&(30,3)&${\tilde a^2_i}$&(22,2)&(25,1)&(20,1.5)\\
\hline
\end{tabular}\\
\end{center}
\begin{center}
Table 5.6 :The demands $\tilde b^p_j$ \\
\begin{tabular}{cccccccccc}
\hline $ j $&1&2&3&4&&1&2&3&4\\
\hline 
 $\tilde b^1_j$ &(10,1.5)&(12,1)&(13,2)&(12,2)& $\tilde b^2_j$ &(5,2)&(5,1.5)&(10,3)&(8,2)\\
\hline
\end{tabular}\\
\end{center}
\begin{center}
Table 5.7 :The transportation capacities $\tilde e_k $\\
\begin{tabular}{cccc}
\hline $k$&1&2\\
\hline 
 $\tilde e_k $&(80,1.5)&(110,2)\\
\hline
\end{tabular}\\
\end{center}
In the following problem, suppose that the confidence levels of constraints are $\gamma^p_i=0.9$;$\beta^p_j=0.9$;$\delta_k=0.9$ and $i = 1,2.3,  j = 1,2,3,4,  k = 1,2, p=1,2$ respectively. From the data tables 5.1, 5.2, 5.3, 5.4, 5.5, 5.6 and 5.7, the corresponding equivalent model to model (6) is constructed as follows:
\begin{equation}
\left\{ \begin{array}{l}
\begin{array}{*{20}{l}}
{E[f_1(\x,\hxi)]=\min \sum\limits_{p= 1}^2 \sum\limits_{i = 1}^3 {\sum\limits_{j = 1}^4 {\sum\limits_{k = 1}^2 {e^{1p}_{ijk}}} {x^p_{ijk}}} }\\
{E[f_2(\x,\hxi)]=\min \sum\limits_{p= 1}^2 \sum\limits_{i = 1}^3 {\sum\limits_{j = 1}^4 {\sum\limits_{k = 1}^2 {e^{2p}_{ijk}}} {x^p_{ijk}}} }
\end{array}\\
s.t.\left\{ \begin{array}{l}
\sum\limits_{j = 1}^4 {\sum\limits_{k = 1}^2 {{x^p_{ijk}} - \left[ {{e^p_i} + \frac{{{\sigma _i}\sqrt 3 }}{\pi }\ln \frac{{1 - {\gamma^p _i}}}{{{\gamma^p _i}}}} \right] \le 0} } ,i = 1,2.3, p=1,2.\\
\left[ {{e'^p_j} + \frac{{{\sigma'^p_j}\sqrt 3 }}{\pi }\ln \frac{{{\beta^p _j}}}{{1 - {\beta^p _j}}}} \right] - \sum\limits_{i = 1}^3 {\sum\limits_{k = 1}^2 {{x^p_{ijk}} \le 0,} } j = 1,2,3,4, p=1,2.\\
\sum\limits_{p = 1}^2\sum\limits_{i = 1}^3 {\sum\limits_{j = 1}^4 {{x^p_{ijk}} - \left[ {{e''_k} + \frac{{{\sigma'' _k}\sqrt 3 }}{\pi }\ln \frac{{1 - {\delta _k}}}{{{\delta _k}}}} \right] \le 0,} } k = 1,2.\\
{x^p_{ijk}} \ge 0, i = 1,2.3., j = 1,2,3,4., k = 1,2., p=1.2.
\end{array} \right.
\end{array} \right.
\end{equation}

In order to obtain a Pareto optimal solutions of model (10), the convex combination method is first used. Thus, the deterministic multi-objective programming problem (10) is transformed as follows based on model (2):  
\begin{equation}
\left\{ \begin{array}{l}
\begin{array}{*{20}{l}}
{\min \left (\lambda_1\sum\limits_{p= 1}^2 \sum\limits_{i = 1}^3 {\sum\limits_{j = 1}^4 {\sum\limits_{k = 1}^2 {e^{1p}_{ijk}}} {x^p_{ijk}}}+\lambda_2 \sum\limits_{p= 1}^2 \sum\limits_{i = 1}^3 {\sum\limits_{j = 1}^4 {\sum\limits_{k = 1}^2 {e^{2p}_{ijk}}} {x^p_{ijk}}} \right) }\\
\end{array}\\
s.t.\left\{ \begin{array}{l}
\sum\limits_{j = 1}^4 {\sum\limits_{k = 1}^2 {{x^p_{ijk}} - \left[ {{e^p_i} + \frac{{{\sigma _i}\sqrt 3 }}{\pi }\ln \frac{{1 - {\gamma^p _i}}}{{{\gamma^p _i}}}} \right] \le 0} } ,i = 1,2.3, p=1,2.\\
\left[ {{e'^p_j} + \frac{{{\sigma'^p_j}\sqrt 3 }}{\pi }\ln \frac{{{\beta^p _j}}}{{1 - {\beta^p _j}}}} \right] - \sum\limits_{i = 1}^3 {\sum\limits_{k = 1}^2 {{x^p_{ijk}} \le 0,} } j = 1,2,3,4, p=1,2.\\
\sum\limits_{p = 1}^2\sum\limits_{i = 1}^3 {\sum\limits_{j = 1}^4 {{x^p_{ijk}} - \left[ {{e''_k} + \frac{{{\sigma''_k}\sqrt 3 }}{\pi }\ln \frac{{1 - {\delta_k}}}{{{\delta_k}}}} \right] \le 0,} } k = 1,2.\\
{x^p_{ijk}} \ge 0, i = 1,2.3., j = 1,2,3,4., k = 1,2., p=1.2.
\end{array} \right.
\end{array} \right.
\end{equation}
where $ \lambda_1  +\lambda_2 =1$ denotes the weight parameters for the objectives.
Then, the  optimal solutions of model (11) with various weight coefficients are obtained. These results are given in table 5.8.
\begin{center}
	Table 5.8 :Results of the convex combination model with some weight coefficients\\
	\begin{tabular}{ccccc} 
		\hline $ \lambda_1$&$\lambda_2$&$E[f_1]$&$ E[f_2] $& optimal value of (11)\\
		\hline 
		1&0&336.964&2232.086&336.964\\
		0.75&0.25&462.541&1725.394&778.255\\
		0.50&0.50&578.579&1531.072&1054.825\\
		0.25&0.75&716.810&1436.487&1256.568\\
		0&1&826.795&1408.9912&1408.991\\
		\hline
	\end{tabular}\\
\end{center}

Following from the compromise programming model (3),  the corresponding programming model to problem (10) is formulated as:
\begin{equation}
\left\{ \begin{array}{l}
\begin{array}{*{20}{l}}
{\min\sqrt{ \left (\sum\limits_{p= 1}^2 \sum\limits_{i = 1}^3 {\sum\limits_{j = 1}^4 {\sum\limits_{k = 1}^2 {e^{1p}_{ijk}}} {x^p_{ijk}}}-E^*_1\right )^2+\left (\sum\limits_{p= 1}^2 \sum\limits_{i = 1}^3 {\sum\limits_{j = 1}^4 {\sum\limits_{k = 1}^2 {e^{2p}_{ijk}}} {x^p_{ijk}}}-E^*_2\right )^2} }\\
\end{array}\\
s.t.\left\{ \begin{array}{l}
\sum\limits_{j = 1}^4 {\sum\limits_{k = 1}^2 {{x^p_{ijk}} - \left[ {{e^p_i} + \frac{{{\sigma _i}\sqrt 3 }}{\pi }\ln \frac{{1 - {\gamma^p _i}}}{{{\gamma^p _i}}}} \right] \le 0} } ,i = 1,2.3, p=1,2.\\
\left[ {{e'^p_j} + \frac{{{\sigma'^p_j}\sqrt 3 }}{\pi }\ln \frac{{{\beta^p _j}}}{{1 - {\beta^p _j}}}} \right] - \sum\limits_{i = 1}^3 {\sum\limits_{k = 1}^2 {{x^p_{ijk}} \le 0,} } j = 1,2,3,4, p=1,2.\\
\sum\limits_{p = 1}^2\sum\limits_{i = 1}^3 {\sum\limits_{j = 1}^4 {{x^p_{ijk}} - \left[ {{e''_k} + \frac{{{\sigma''_k}\sqrt 3 }}{\pi }\ln \frac{{1 - {\delta_k}}}{{{\delta_k}}}} \right] \le 0,} } k = 1,2.\\
{x^p_{ijk}} \ge 0, i = 1,2.3., j = 1,2,3,4., k = 1,2., p=1.2.
\end{array} \right.
\end{array} \right.
\end{equation}
where$  E^*_1 $ and $ E^*_2 $ are the individual optimal values of the objective functions in problem (10) ignoring the other one. According to Table $ 5.8 $, it is easy to determine (minimum of the objectives ) that $ E^*_1=336.964$ and $ E^*_2=1408.991.$

By using the Maple 18.02 optimization toolbox, the compromise programming model (12) is solved and then results are summarized in table 5.9.    
\begin{center}
	Table 5.9 :Optimal solution of problem  (12) with minimizing distance function method\\
	\begin{tabular}{ccccc} 
		\hline optimal solutions&$E[f_1]$&$E[f_2] $&objective value of model (12) \\
		\hline 
		$ x^1_{121} = 10.319,x^1_{122} = 2.893, x^1_{232} = 0.874, x^1_{241} = 14.423, $\\ $ x^1_{311} = 11.817, x^1_{331}= 14.549, x^2_{112} = 2.337, x^2_{121} = 3.817, $\\ $ x^2_{122} = 3, x^2_{142} = 10.423, x^2_{232} = 13.634, $\\ $ x^2_{311} = 2.650, x^2_{312} = 2.435 $&551.148&1571.781&257.418\\
		\hline
	\end{tabular}\\
\end{center}

Consequently, the Pareto optimal solution  of model (10) can be selected based on results of Table 5.8 and 5.9. In addition, the minimizin distance function method in this paper is generated the best solution according to the convex combination method.

\section{Conclusions}

In this paper, a multi-objective multi-item solid transportation problem with uncertain parameters is investigated. Using the advantage of uncertainty theory, the model is first transformed into the deterministic multi-objective multi-item solid transportation problem. Then the deterministic model based on the expected value of each objective under the chance constraints is reduced to single objective programming problems by employing convex combination method and minimizing distance function method. The solution procedure for each method is illustrated with a numerical example and the results derived from both methods were compared. Finally, Pareto optimal solution of the considered model (10) is obtained.    

It is noted that uncertain programming can be studied so extensively in actual applications since multi-objective and uncertain criteria broadly arise in all types of real-world optimization problems.

\end{document}